%From fuchino@leibniz.math.fu-berlin.de Sun May  7 18:15:59 1995
%Received: from ki1.Chemie.FU-Berlin.DE by sunset.ma.huji.ac.il with SMTP id AA16170
%  (5.67b/HUJI 4.152 for <shlhetal@math.huji.ac.il>); Sun, 7 May 1995 18:09:08 +0300
%Received: by ki1.chemie.fu-berlin.de (Smail3.1.28.1)
%	  from leibniz.math.fu-berlin.de (160.45.40.10) with smtp
%	  id <m0s87vv-0000YnC>; Sun, 7 May 95 17:07 MEST
%Received: by leibniz.math.fu-berlin.de (/\=-/\ Smail3.1.18.1 #18.14)
%	  id <m0s87vu-000sLvC>; Sun, 7 May 95 17:07 MET DST
%Message-Id: <m0s87vu-000sLvC@leibniz.math.fu-berlin.de>
%Date: Sun, 7 May 95 17:07 MET DST
%From: fuchino@math.fu-berlin.de (Sakae Fuchino)
%To: shlhetal@math.huji.ac.il
%In-Reply-To: <199505061459.AA29899@sunset.ma.huji.ac.il> (message from Saharon Shelah's office on Sat, 6 May 1995 17:59:23 +0300 (EET DST))
%Subject: Re: joint paper with Saharon
%Status: OR
%
%Dear Andrzej,
%the following is the version which we think is final of 556:
%
%-- 
%     Best Regards,
%                    Saka\'e     $@!!!J^<Ln(J $@>;!K(J
%
%     fuchino@math.fu-berlin.de
%     ( fuchino@logic.info.waseda.ac.jp )
%--------- cut here -----------
%%% Return-Path: <sabina@somerville.math.fu-berlin.de>
%%% From: sabina@math.fu-berlin.de (Prof. Sabine Koppelberg)
%%% Subject: paper
%%% To: fuchino@math.fu-berlin.de (Sakae Fuchino)
%%% Date: Fri, 5 May 1995 14:38:46 +0200 (MET DST)
%%% X-Mailer: ELM [version 2.4 PL22]
%%% MIME-Version: 1.0
%%% Content-Type: text/plain; charset=US-ASCII
%%% Content-Transfer-Encoding: 8bit
%%% Content-Length: 28066     

\input amstex
\documentstyle{amsppt}
\ifx\kanjiskip\otankonasu\magnification \magstep1\fi %
%%%
%%% coded by S.K.
%%% extended by S.F. on: 280495
%%% extended further by S.K. on: 150595
\topmatter
%%%%%%%%%%%%%%%%%%%%%%%%%%%%%%%%%%%%%%%%%%%%%%%%%%%%%%%%%%%%%%%
\title \nofrills   A game on partial orderings        \endtitle              
\author  Saka\'e Fuchino, Sabine Koppelberg, 
and Saharon Shelah
\endauthor
%%%%%%%%%%%%%%%%%%%%%%%%%%%%%%%%%%%%%%%%%%%%%%%%%%%%%%%%%%%%%%%
%\affil            \endaffil
\address
{\it Saka\'e Fuchino, Sabine Koppelberg}\newline
Institut f\"ur Mathematik II, Freie Universit\"at Berlin, Arnimallee 3, 
14195 Berlin, Germany.\medskip
\noindent{\it Saharon Shelah}\newline
Institute of Mathematics, The Hebrew University of Jerusalem, 91904 
Jerusalem, Israel.\newline
Department of Mathematics, Rutgers University, New Brunswick, NJ 08854, U.S.A.
\endaddress
%\email           \endemail
\dedicatory      \enddedicatory
\date May 6, 1995             \enddate
\thanks  This research was begun when the first author was working at the 
Hebrew University of 
Jerusalem. He would like to thank The Israel Academy of Science and 
Humanities for enabling his stay there. 
The second author was partially supported by the Deutsche 
Forschungsgemeinschaft (DFG) Grant No.\ Ko 490/7--1. 
The third author would like to thank ``Basic Research Foundation'' of The 
Israel Academy of Science and Humanities and the Edmund Landau Center for 
research in Mathematical Analysis supported by the Minerva Foundation 
(Germany) for their support. 
\endthanks 
%\translator          \endtranslator
\keywords  Freese-Nation property, games, trees, linear orders, Boolean 
algebras 
\endkeywords 
\subjclass 03E05   \endsubjclass
%-------------------------------------------------Macros added by S.F.--------
\def\fnsp#1#2{\text{}^{#1^{\text{}\!}}#2}
\def\calH{{\Cal H}}
\def\cardof#1{{\mathopen{\mid}{#1}\mathclose{\mid}}}
\def\setof#1#2{\{#1:#2\}}
\def\mapping#1#2#3{#1:#2\rightarrow #3}
\def\cf{\text{\rm cf\,}}
\def\ran{\text{\rm range\,}}
\def\st{such that}
%-----------------------------------------------------------------------------
\abstract          
We study the determinacy of the game $G_\kappa(A)$ introduced in 
\cite{FuKoShe} for uncountable regular $\kappa$ and several classes of 
partial orderings $A$. Among trees or Boolean algebras, we can always find 
an $A$ \st\ $G_\kappa(A)$ is undetermined. For the class of linear orders, 
the existence of such $A$ depends on the size of $\kappa^{<\kappa}$. In 
particular we obtain a characterization of $\kappa^{<\kappa}=\kappa$ in 
terms of determinacy of the game $G_\kappa(L)$ for linear orders $L$. 
\endabstract

\endtopmatter

\document

We consider in this paper the question whether for every partially ordered set 
$(A, \leq )$, the game  $G_\kappa(A)$  described below is 
determined, i.e.\ whether one of the players 
has a winning strategy. Here and in the following, except for the motivation given
below,  $\kappa$ is always a regular uncountable  cardinal.  More precisely we study
the question for trees,   Boolean algebras and linear orderings. In fact there are
trees, resp.   Boolean algebras, $A$ of size $\kappa^+$ for which $G_\kappa(A)$ is
not  determined (Propositions 6 and 11); for linear orders, the situation is  more
complex:  if $\kappa^{<\kappa}=\kappa$, then for every linear order $L$, 
$G_\kappa(L)$ is determined (Proposition 2); otherwise there is a linear 
order $L$ of size $\kappa^+$ such that $G_\kappa(L)$ is not determined 
(Proposition 8).

The motivation for this question comes from the paper \cite{FuKoShe} which
 in turn was motivated by  \cite{HeSha}. A Boolean algebra  $A$  is said to have the
 Freese-Nation property if there exists a function  $f$  which assigns to every
 $a \in A$  a finite subset   $f(a)$  of   $A$  such that if  $a, b \in A$  
satisfy  $a \leq b$, then  $a \leq x \leq b$  holds for some  
$x \in f(a) \cap f(b)$. This property is closely related to projectivity; 
in fact, every projective Booleran algebra has the Freese-Nation property
 (but not conversely). Heindorf proved that the Freese-Nation property is
equivalent to open-generatedness, a notion originally introduced in topology by
\v S\v cepin. In \cite{FuKoShe}, it is generalized to from  $\omega$  to 
regular 
cardinals  $\kappa$  and from Boolean algebras to arbitrary partial 
orderings. This generalization is called $\kappa$-Freese-Nation property and 
the following equivalence was proved: a partial ordering 
$A$ has the  $\kappa$-Freese-Nation property iff there is a  closed unbounded
subset  $\Bbb C$  of  $[A]^\kappa$  such that  $C \leq_\kappa A$  holds 
for all 
$C \in \Bbb C$  iff  in the game  $G_\kappa(A)$,  Player II has a winning 
strategy. In fact, in all examples considered in
\cite{FuKoShe}, either I or II  has a winning strategy.
 
Let us define the game  $G_\kappa(A)$  and some relevant notions for a partial
ordering  $A$.  $X \subseteq A$  is said to be cofinal  (coinitial) in 
$A$  if, for every  $a \in A$,  there is some $x \in X$  such that 
$a \leq x$  ($a \geq x$). 
$\text{cf } A$  resp.\  $\text{ci } A$   is the smallest cardinality of a 
cofinal 
resp.\ coinitial subset of  $A$.

For  $R \subseteq A$  and  $a \in A$, we write  $R \uparrow a$ for the set 
$\{x \in R : a \leq x \}$  and 
$R \downarrow a$  for  $\{x \in R : x \leq a \}$. The {\it type} of  $a$ 
over  $R$ 
is the pair

$$\text{tp }(a,R) = (\text{cf } R \downarrow a ,\text{ci } R \uparrow a  ). $$

$R \subseteq A$ is said to be a 
$\kappa$-subset or a  $\kappa$-substructure of  $A$, written 
$R \leq_{\kappa} A$, if for all  $a \in A$, the sets 
$R \downarrow a$  and  $R \uparrow a$
have cofinality resp.\ coinitiality less than  $\kappa$.

The game   $G_\kappa(A)$  is played on  $A$  as follows. Players  I  
and II  alternatively choose an increasing chain of 
subsets $x_\alpha$  and  $y_\alpha$  of 
$A$  for  $\alpha < \kappa$  (i.e.\ I  chooses  $x_0$,   
II  chooses  $y_0$, I  chooses  $x_1$, II  chooses  $y_1$, etc.) such that  
$x_\alpha$  and  $y_\alpha$  have size less than  $\kappa$,  
$x_\alpha \subseteq y_\alpha$  and  $\bigcup_{\nu <\alpha} y_\nu 
\subseteq x_\alpha$. 
In the end of a play,
II  wins iff the result   $R = \bigcup_{\alpha <\kappa } x_\alpha = 
\bigcup_{\alpha <\kappa } y_\alpha$  of the play
is a  $\kappa$-subset of  $A$.

Note that in this game, Player  II  has a winning strategy for any partial 
ordering   $A$  of size at most  $\kappa$: she can play so that every element
of  $A$  is gradually captured in one of the  ${y_\alpha} \rq$s.

The main body of the paper is organized as follows. In 5., we define a tree  $T = T(S)$,
depending on a subset   $S$  of  $\lambda = \kappa^+$. If neither  $S$  nor  $\lambda
\setminus S$  are in the ideal  $I_\lambda$  defined in  3., then  $T$  is not determined
(Proposition 6). From  $T$, we define a linear order  $L_T$ in 7. and a Boolean algebra 
$B_T$  in 10. such that  $G_\kappa (L_T)$  and  $G_\kappa (B_T)$  are not determined
(Propositions  8 and 11). The construction of  $L_T$  requires the extra assumption 
 $\kappa^{<\kappa} > \kappa$  --- cf.\ Proposition 2.

Let us start with an easy example. 

\example{1. Example} If  $\kappa^+$ (or $(\kappa^+)^{-1}$, the reverse 
order type of $\kappa^+$) 
embeds into  $A$, then Player I has a winning
strategy in  $G_\kappa (A)$: assume, for simplicity, that 
$\kappa^+\subseteq A$. 
We define a partial function  $f$  from  $A$  into  $\kappa^+$  by letting 
$f(a)$  for  $a \in A$ be the least  $\alpha \in \kappa^+$  such that 
$a \leq \alpha$, if such an  $\alpha$  exists. Clearly  $f$  is order 
preserving 
and satisfies $f(a) = a$  for  $a \in \kappa^+$. Player I wins by assuring 
that 
the result  $R$ of a play satisfies

(a) $R \cap \kappa^+$  has cofinality  $\kappa$

(b)  if  $a \in R$  and  $f(a)$  exists, then  $f(a) \in R$.  
\qed
\endexample 

The following proposition shows that the assumption 
$\kappa^{<\kappa}>\kappa$ in 6.\ and 7.\ cannot be dispensed with. 

\proclaim{2. Proposition} 
Assume that $\kappa^{<\kappa}=\kappa$. If $(L,<_L)$ is a linear 
order of cardinality $>\kappa$, then Player I has a winning strategy in 
$G_\kappa(L)$. 
Hence the game $G_\kappa(L)$ is 
determined for any linear order $L$ under $\kappa^{<\kappa}=\kappa$. 
\endproclaim
\demo{Proof} 
Let $\chi$ be sufficiently large.  $\calH(\chi)$  denotes the set of all sets  
which are hereditarily of size less than  $\chi$. 
We show:
\proclaim{Claim}
Suppose $M$ is an elementary submodel of $(\calH(\chi),  \in)$  
such that $(L,<_L)\in M$ 
and $\fnsp{\kappa>}{M}\subseteq M$. Then for any 
$d \in L\setminus M$, either $L\cap M\downarrow d $ has cofinality 
$\geq\kappa$ or $L\cap M\uparrow d $ has coinitiality $\geq\kappa$. 
\endproclaim
\demo{Proof of the Claim}
Otherwise, some 
$d\in L\setminus M$ fills a gap $(X,Y)$ in $L\cap M$ \st\ $\cardof{X}$, 
$\cardof{Y}<\kappa$ and $(X,Y)$ is unfilled  
inside $L\cap M$. But $(X,Y)\in M$ by 
$\fnsp{\kappa>}{M}\subseteq M$ 
and $M\prec\calH(\chi)$, a contradiction.
\qed
\enddemo
Now Player I wins in $G_\kappa(L)$ by choosing an increasing sequence 
$M_\alpha$, $\alpha<\kappa$, of elementary submodels of $\calH(\chi)$ 
along with his moves $x_\alpha$, $\alpha<\kappa$, such that 
$(L,<_L)\in M_0$, 
$x_\alpha\subseteq M_\alpha$, 
$\fnsp{\kappa>}{M_\alpha}\subseteq {M_\alpha}$, 
$\cardof{M_\alpha}=\kappa$ and 
$\bigcup_{\alpha<\kappa}x_{\alpha}=M\cap L$ where 
$M=\bigcup_{\alpha<\kappa}M_\alpha$. 
Such a  choice is possible because of our assumption 
$\kappa^{<\kappa}=\kappa$. The result of the game 
$L\cap M$ is not a $\kappa$-subset of $L$, by the Claim above.
\qed        
\enddemo

\example{3. The ideal $I_\lambda$}
For the rest of the paper, fix $\lambda=\kappa^+$ (where $\kappa$ was 
a regular uncountable cardinal). Let us first recall the definition and some 
properties of the ideal $I_\lambda$ on $\lambda$ introduced by Shelah, see 
e.g.\ \cite{She,Chapter VIII}. 
Fix a 
sufficiently large cardinal $\chi>\lambda$; we work in the structure 
$(\calH(\chi),\in,<^*)$ where $<^*$ is some fixed well-ordering of 
$\calH(\chi)$. 
For $x\in\calH(\chi)$ and $\gamma<\lambda$, call $(M_i)_{i<\kappa}$ an 
$x$-approximation of $\gamma$ if:
\roster
\item $M_i\prec (\calH(\chi),\in,<^*)$, $\cardof{M_i}<\kappa$
\item $x,\lambda\in M_0$
\item $(M_i)_{i<\kappa}$ is a continuously increasing chain 
\item $(M_i)_{i\leq j}\in M_{j+1}$ for all $j<\kappa$
\item $M=\bigcup_{i<\kappa}M_i$ satisfies $M\cap\lambda=\gamma$.
\endroster
For $x\in\calH(\chi)$, put
$C_x=\setof{\gamma\in\lambda}{\text{there is an }x\text{-approximation of }
\gamma}$
and define $I_\lambda$ by
$$
I_\lambda=\setof{A\subseteq\lambda}{A\cap C_x=\emptyset,\text{ for some }
x\in\calH(\chi)}.
$$
It is not difficult to check that $I_\lambda$ is a $\lambda$-complete 
proper ideal containing all singletons and that 
$$
N=\setof{\gamma\in\lambda}{\cf\gamma=\kappa}\in {I_\lambda}^*
$$
(i.e.\ $\lambda\setminus N\in I_\lambda$). By Ulam's Theorem (cf.\ 
\cite{Je, 27.8}), every  $A\subseteq\lambda$ not in $I_\lambda$ can be represented 
as the disjoint union $A=A_1\cup A_2$  where  $A_1$, $A_2\not\in I_\lambda$. 
\endexample 

\example{4. The game $G_\kappa(T)$ for a tree $T$}
Assume that $(T,<_T)$ is a tree of height $\kappa+1$. We call 
$Y\subseteq T$ a {\it subtree of $T$} if for all $y\in Y$ and $x<_T y$, also 
$x\in Y$. $Y$ is {\it closed in $T$} if the following holds: if $x\in T$ 
is in the $\kappa$'th level and all predecessors of $x$ are in $Y$, then 
$x\in Y$.

In $G_\kappa(T)$ each of the players can ensure that the result $Y$ of a 
play will be a subtree of $T$. And in this case, Player II wins, i.e.\ 
$Y\leq_\kappa T$, iff $Y$ is closed in $T$. 
\endexample

\example{5. Construction of the tree $T=T(S)$}
Recall that $\lambda=\kappa^+$ and 
$N=\setof{\gamma\in\lambda}{\cf\gamma=\kappa}$. 
Depending on a subset $S$ of $N$, we construct a tree $T=T(S)$; in 
fact, we shall show that if  $T=T(S)$  where  $S\subseteq N$ and  $S$, 
$N\setminus S\not\in I_\lambda$, then none of the players has a winning strategy.

Assume $S\subseteq N$. For each $\gamma\in S$, fix a function 
$\mapping{f_\gamma}{\kappa}{\gamma}$ \st\ $\ran f_\gamma$ is cofinal in 
$\gamma$. Let 
$$
T=T(S)=\setof{f_\gamma \restriction\alpha}{\gamma\in S,\,\alpha\leq\kappa},
$$
a tree under set-theoretic inclusion. Clearly $T$ has height $\kappa+1$ 
if $S$ is nonempty, $\setof{f_\gamma}{\gamma\in S}$ is the $\kappa$'th 
level of $T$, and $\cardof{T}=\lambda$ if $\cardof{S}=\lambda$.
\endexample

\proclaim{6. Proposition}
Let $T=T(S)$ for $S\subseteq N$.\medskip

(a) If $S\not\in I_\lambda$, then Player II has no winning strategy in 
$G_\kappa(T)$. 

(b) If $N\setminus S\not\in I_\lambda$, then Player I has no winning 
strategy in $G_\kappa(T)$. \medskip

\noindent
Thus if both $S$ and $N\setminus S$ are not in $I_\lambda$, then 
the game $G_\kappa(T)$ is undetermined. 
\endproclaim
\demo{Proof}
(a) Suppose that $\sigma$ is a strategy for Player II; we show that it  is not a
winning strategy. Let  $x=(\sigma,(f_\gamma)_{\gamma \in S})$. Since 
$ S\not\in I_\lambda$, 
there is a $\delta \in  S \cap C_x$; let $(M_i)_{i<\kappa}$ be an 
$x$-approximation of $\delta$. In a game in which Player II plays 
according to $\sigma$, Player I can 
ensure that the result  $Y\subseteq T$  of the play will be the subtree
$$
Y=\setof{f_\gamma\restriction\alpha}{\gamma\in S\cap\delta,\,\alpha \leq \kappa}.
$$
More precisely, in the  $i$'th move, Player I may take a subset  $x_i$  of 
$T\cap M_{i+1}$ so that all 
elements of $Y$ are gradually captured. Furthermore, using the well-ordering  $<^*$,
Player I can  ensure that each of his moves $x_i$ is definable so that 
$(x_j,y_k)_{j\leq i,\,k<i}$ and hence also the next      move 
$\sigma((x_j,y_k)_{j\leq i,\,k<i})$ by Player II will be an element 
of $M_{i+1}$. 

Now  $\delta \in S$  and thus  $f_\delta$   witnesses that  $Y$  is not closed in 
$T$, i.e. Player  I  wins.

The proof of (b) is similar to (a). If Player I plays according to a 
strategy  $\tau$, Player II can assure that the result 
$Y\subseteq T$  has the form
$Y=\setof{f_\gamma\restriction\alpha}{\gamma\in S\cap\delta,\,\alpha \leq \kappa}$ 
for some $\delta\in N \setminus S$. Thus  $Y$  is closed in  $T$  and Player II 
wins.        \qed
\enddemo

\example{7. Construction of the linear order $L_T$} Assume that 
$\kappa^{<\kappa}>\kappa$; let  $(T, <_T)$  be 
any tree of height
$\kappa + 1$  and size  $\lambda=\kappa^+$. We shall construct a linear order
$L=L_T$  of size $\lambda$. Moreover, we shall define for  $Y \subseteq T$ a 
subset  $L_Y$ of  
$L$  such that  $\cardof{L_Y}=\cardof{Y}$ holds for infinite  $Y$ \st,
in the game  $G_\kappa(L)$, 
each player can ensure that the result  $R$  has the form  $L_Y$  for  $Y$  a 
subtree of  $T$.

Let us first note that there exists a linear order   $I$  of size  $\lambda$ 
without any sequences (i.e.\ increasing or decreasing sequences) 
of type $\kappa$. 
This holds because our assumption $\kappa^{<\kappa}\geq\kappa^+=\lambda$
implies that  
$\lambda \leq 2^\mu $, for some  $\mu < \kappa$, and the lexicographic 
ordering 
on  ${}^\mu 2$  has no sequence of type  $\mu^+$ (cf. \cite{Je, 29.4}), hence 
no 
sequence of type  $\kappa$. It follows that, letting $I$ be any 
subordering of $\fnsp{\mu}{2}$ of cardinality $\lambda$, every subset 
of $I$ has cofinality and coinitiality less than $\kappa$. 

The following notation concerning the tree  $(T, <_T)$ will be used in the rest 
of 7.\ and in 8.: for  $\alpha \leq \kappa$, $\text{lev}_\alpha T$ is 
the 
$\alpha$\rq th level of $T$. For   $t \in T$,  $\text{pred } t$   
is the set of predecessors of  $t$  in 
$T$  and  $\text{ht } t$ is the height of  $t$. 
%% baka
For $\alpha \leq \text{ht } t $, 
$\text{pr}_\alpha t$, the projection of  $t$  to level   $\alpha$,  
is the unique predecessor of  $t$  in the 
$\alpha$\rq th level. Call  $x$, $y \in T$  equivalent and write 
$ x \sim y$  if  $\text{pred } x  = \text{pred } y $  
and let  $\overline x$  be the equivalence class of  $x$. For each 
equivalence class $\overline{x}$, 
since  $\vert \overline x \vert \leq \lambda$, we can fix a linear order  
$\leq_{\overline x}$  on  $\overline x$  without any sequences of type  
$\kappa$.

The linear order we construct is a sort of squashing of $T$ with respect 
to $\leq_{\overline{x}}$, $x\in T$: 
we put  $L = \{a_t, b_t : t \in T  \}$  where the elements  $a_t$, 
$b_t$,  $t \in T$, are all pairwise distinct.
The linear order  $<_L$  on  $L$  is defined as follows: 
we will have  $a_t <_L b_t$  for all  $t \in T$. 
Now assume  $x$,  $y \in T$. If  $x <_T y$, then we put 
$a_x <_L a_y <_L b_y <_L b_x$. 
If  $x$  and  $y$  are incomparable in  $T$, 
let  $\alpha \leq \kappa$  be minimal such that  $\text{pr}_\alpha x \neq 
\text{pr}_\alpha y$; thus  
$\text{pr}_\alpha x \sim \text{pr}_\alpha y$. Then if 
$\text{pr}_\alpha x  <_{\overline{\text{pr}_\alpha x }} \text{pr}_\alpha y$,
we let $a_x <_L b_x <_L a_y <_L b_y$. Finally, for $Y\subseteq T$ let 
$L_Y=\setof{a_t,b_t}{t\in Y}$. 
\endexample

\proclaim{8. Proposition} If  $Y$  is a subtree of  $T$, then  
$L_Y \leq_\kappa L_T$  
iff  $Y$  is closed in  $T$.  
In particular, if $G_\kappa(T)$ is undetermined, then so is 
$G_\kappa(L_T)$. 
\endproclaim

>From Propositions 2, 6, and 8 (plus the observation in 7.\ that
$\kappa^{<\kappa}>\kappa$  implies the existence of a linear order of size 
$\lambda$  without sequences of type $\kappa$), we obtain
the  following equivalences to the condition  $\kappa^{<\kappa}=\kappa$.

\proclaim{9. Corollary} Let  $\kappa$  be a regular uncountable cardinal.  \par
(a)  If  $\kappa^{<\kappa}>\kappa$, then there is  a linear order $L$ of cardinality
$\lambda=\kappa^+$ \st\ $G_\kappa(L)$  is undetermined.   \par
(b)  The following are equivalent:
\roster
\item $\kappa^{<\kappa}=\kappa$;
\item in every linear order of cardinality $>\kappa$, there is an increasing or a
decreasing sequence of order type $\kappa$;
\item $G_\kappa(L)$ is determined for every linear order $L$.\qed
\endroster
\endproclaim

Let us explain how the second assertion of Proposition 8 follows from
 the first one: each of the players in  $G_\kappa(L_T)$  (say II, playing against 
some strategy  $\tau$  of  Player I) can ensure that the result of the play 
is  $R = L_Y$, for some subtree  $Y$  of  $T$. Playing simultaneously on  $T$ as in
the proof of Proposition 6, she can ensure that  $Y$  is closed in  $T$. Thus 
$R = L_Y$  is a  $\kappa$-substructure of  $L_T$  and  II  wins. The same reasoning
applies, of course, to the proof of Proposition 11.

\demo{Proof of Proposition 8}  Suppose first that  $Y$  is not closed and 
pick some  $t$  in the
highest level  $K$  of  $T$  such that  $t \notin Y$  but  $\text{pred } t 
\subseteq Y$. 
Then  $\{a_y: y \in \text{pred } t \}$  is an increasing sequence of type  
$\kappa$,
and it is a cofinal subset of  $L_Y \downarrow l$  where  $l = a_t$. Thus  
$L_Y$  is
not a  $\kappa$-subset of  $L$.

Now assume that  $Y$  is closed in  $T$  and fix  $l \in L \setminus L_Y$. 
We have
to analyze the cofinality of  $L_Y \downarrow l$  and the coinitiality of  
$L_Y \uparrow l$; by symmetry, we will consider 
$ \text{cf } (L_Y \downarrow l)$. 
Now let  $l = a_t$  or  $l = b_t$  for some  $t \in T \setminus Y$; since  $Y$  
is a subtree of  $T$,  $a_t$  and  $b_t$  realize the same cut in  $L_Y$. 
Thus we assume that  $l = a_t$.

We may also assume that  $\text{ht } t < \kappa$  and 
$\text{pred } t \subseteq Y$. 
For this, consider the least  element  $t^*$  of 
$\text{pred } t \setminus Y$. 
Now  $\text{ht } t^* < \kappa$  since  $Y$  is a closed subtree of  $T$; 
moreover,  
$a_t$  and  $a_{t^*}$  realize the same cut in  $L_Y$. Thus we consider 
$t^*$ instead of  $t$.

To prove  $\text{cf } (L_Y \downarrow l) < \kappa$, consider the following 
subsets
of  $L$ respectively  $Y$: let

$$N = \{a_x: x <_T  t \} ;$$
thus  $N$  is a subset of  $L_Y \downarrow l$  of size  less than  $\kappa$. 
Next, put  $\gamma = \text{ht } t$  and

$$ Y^\prime = \{z \in Y: z \in \text{lev}_\gamma T, z \sim t, z <_{\overline 
t} t \}.$$
$Y^\prime$  is included in the  $\sim$-equivalence class of  $t$, thus it 
has a cofinal subset  $Y^{\prime \prime}$  of size  less than  $\kappa$. We 
put

$$ N^\prime = \{ b_z : z \in Y^{\prime \prime}\},$$
again a  subset of   $L_Y \downarrow l$  of size  less than  $\kappa$.

We prove that  $N \cup N^\prime$  is cofinal in  $L_Y \downarrow l$. For, 
let 
$x \in L_Y$  and  $x <_L l$,
say  $x \in \{a_y, b_y\}$  where  $y \in Y$. Consider the relative position 
of $t$  and  $y$  in  $T$. It is impossible that  $t <_Ty$, since  $Y$  is 
a subtree of  $T$  and  $t \notin Y$. 

If  $y <_T t$, then  $a_y <_L a_t = l <_L b_t <_L b_y$  holds, hence 
$x = a_y \in N$. Otherwise, let  $\alpha$  be minimal such that  
$\text{pr}_\alpha y \neq \text{pr}_\alpha t$;
thus  $\alpha \leq \gamma$.

If  $\alpha < \gamma$, then let  $z = \text{pr}_\alpha  t$; it follows that 
$x \leq_L b_y <_L a_z \in N$. Otherwise  $\alpha = \gamma$,  
$ \text{pr}_\alpha y \sim t $  and hence  $ \text{pr}_\alpha y \in 
Y^\prime$. Take
 $z \in Y^{\prime \prime}$  such that 
$ \text{pr}_\alpha y \leq_{\overline t}  z$;  then 
$x \leq b_y \leq b_z \in N^\prime$.\qed 
\enddemo

\example{10. Construction of the Boolean algebra $B_T$}   Let 
$(T, <_T)$  be any 
tree of height
$\kappa + 1$  and size  $\lambda$. We shall construct a Boolean algebra
$B_T$  of size $\lambda$. Moreover, we shall define for  $Y \subseteq T$ a 
subalgebra  $B_Y$ of  
$B_T$  such that  $\vert B_Y  \vert = \vert Y \vert$ holds for infinite  $Y$. 
In the game  $G_\kappa(B_T)$, 
each player can ensure that the result  $R$  has the form  $B_Y$  for  $Y$  a 
subtree of  $T$.       

In fact, we define  $B_T$  to be the Boolean algebra generated by a set  
$\{ x_t: t \in T \}$  freely 
except that  $s \leq_T t$  implies  $x_s \leq x_t$. More precisely, let  
$\text{Fr } (x_t: t \in T)$  be the free Boolean algebra over 
$\{x_t: t \in T\}$, 
let  $B_T$  be the quotient algebra  $\text{Fr } (x_t: t \in T) / K$  
where  $K$  is the ideal of  $\text{Fr } (x_t: t \in T)$  generated by  
$\{x_s \cdot -x_t : s \leq_T t\}$  and let 
$\pi : \text{Fr } (x_t: t \in T) \rightarrow B_T$  be 
the canonical homomorphism. We write  $x_t$  ($\in B_T$)  for  $\pi(x_t)$, 
since  $\pi$ is one-one on the generators  $x_t$  (see the proof of 10.\ 
below). 
For  $Y \subseteq T$, we define  $B_Y$  to be the subalgebra of  $B_T$  
generated by  
$\{x_t: t \in Y\}$.
\endexample

\proclaim{11. Proposition}  If   $Y$  is a subtree of  $T$, then  $B_Y 
\leq_\kappa B_T$ iff  $Y$  is closed in  $T$.  In particular, if 
$G_\kappa(T)$ is undetermined, then so is $G_\kappa(B_T)$. 
\endproclaim

\demo{Proof}   We start with a normal form lemma on the generators of  
$B_T$.

{\it Step 1.}  Let  $w \subseteq T$  be finite and assume  $f: w \rightarrow 
2$. 
Then the elementary product  $q_f = \prod_{f(t) = 1} x_t \cdot \prod_{f(t) 
= 0} -x_t$ 
is nonzero in  $B_T$ iff   $f$  is monotone, i.e.\ $s \leq_T t$  in  $w$  
implies 
$f(s) \leq f(t)$. --- This follows immediately from the definition of the 
ideal $K$  of  $\text{Fr } (x_t: t \in T)$  in 9.

{\it Step 2.}  If  $Y \subseteq T$  is not closed, then  $B_Y$  is not a  
$\kappa$-
subalgebra of  $B_T$.

To see this, fix an element  $t$  in the highest (i.e.\ $\kappa$\rq th) 
level  of $T$
such that  $t \notin Y$  but all predecessors of  $t$  in  $T$  are in  $Y$  and 
consider the
ideal  $I = B_Y \downarrow x_t$  of  $B_Y$. The set  $J = \{ x_s:s <_T t \}$  
is a chain of order type  $\kappa$  included in  $I$; we show that  $J$  
generates  $I$ as an ideal. Thus suppose  $x \in I$  with the aim of 
finding some  $s <_T t$ such 
that  $x \leq x_s$. We may assume that  $x$  is a non-zero elementary 
product  $q_f$  
where  $f: w \rightarrow 2$. By  $q_f \leq x_t$  and Step 1, it follows that  
$f$  
is monotone but  $f \cup \{ (t,0)\}$  is not. Hence there is some  $s \in w$  
such that  $s <_T t$  and  $ f(s) = 1$; thus  $x = q_f \leq x_s$.

{\it Step 3.}  The following remark simplifies Step 4: assume  $B$  is a 
Boolean algebra, 
 $A$  a subalgebra and  $M$,  $N$  are finite subsets of  $B$  such that for 
all
 $m \in M$  and  $n \in N$, there is an element  $\alpha$  of  $A$  
separating  $m$  
and  $n$, i.e.\ we have  $m \leq \alpha$  and  $n \leq -\alpha$  or  $n \leq 
\alpha$  
and  $m \leq -\alpha$. Then there is an  $a \in A$  separating  $\sum M$  
and  $\sum N$:
simply let  $a = \prod_{n \in N} \sum_{m \in M} a_{mn}$  where  $a_{mn} 
\in A$ is such
that   $m \leq a_{mn}$  and  $n \leq -a_{mn}$.

{\it Step 4.}  If  $Y$  is a closed subtree of  $T$, then  $B_Y \leq_\kappa 
B_T$. 

For the proof, fix an element  $b$  of  $B_T$  and consider the ideal
$$ I = \{x \in B_Y  :x \cdot b = 0 \} $$
of  $B_Y$. We shall find  $Z \subseteq T$  such that  $\vert  Z \vert < 
\kappa$  and 
each element of  $I$  is separated from  $b$  by an element of  $B_Z$; since  
$\vert B_Z \vert < \kappa$, this shows that  $I$  is generated by less than  
$\kappa$  
elements.

Fix a finite subset of  $T$  generating  $b$, say
$$b \in \langle x_{s_1}, \dots,x_{s_n}, x_{t_1}, \dots, x_{t_m}  \rangle $$
where every  $s_i$  is in  $Y$ and every  $t_j$  is in  $T \setminus Y$. We 
put
$$ Z = \{s_1, \dots, s_n \} \cup \bigcup \{\text{pred } t_j \cap Y: 1 \leq j 
\le m\} $$
where, for  $t \in T$, $\text{pred } t$  is the set of predecessors of  $t$  
in the tree  $(T, <_T)$.  $Z$ has size less than  $\kappa$  since  $Y$  is closed 
and 
a subtree of  $T$.

Now let  $x \in I$  with the aim of finding an element of  $B_Z$  which 
separates  $x$  
and  $b$. By Step 3, we may assume that both $b$  and  $x$  are 
elementary products 
over the generators of  $B_T$, say 
$$ b = q_h  \text{, }  h: \{s_1, \dots,s_n, t_1, \dots, t_m\} \rightarrow 2  $$
$$ x = q_f \text{, } f:w \rightarrow 2 \text{, } w \subseteq Y $$
where  $h$  and  $f$  are monotone. Define
$$ h^\prime = h \restriction \{s_1, \dots,s_n \} \text{, } 
                                    f^\prime = f \restriction (w \cap Z);$$
we show that either  $q_{h^\prime}$  or  $q_{f^\prime}$  separate  $x$ and  
$b$.

{\it Case 1.}  $f \cup h^\prime$ is not a function or not monotone. --- Then  
$b \leq q_{h^\prime}$  and  $x \cdot q_{h^\prime} = 0$. 

Note that if Case 1 does not hold, then also  $f \cup h$  is a function: 
otherwise, 
let  $r \in w \cap \{s_1, \dots,s_n, t_1, \dots, t_m\}$  be such that  $f(r) 
\neq h(r)$.
 Then  $r \in Y$  and thus  $r = s_i$  for some  $i$, hence  $r \in \text{dom 
} f 
 \cap \text{dom } h^\prime $. Note also that, since  $x \cdot b = 0$,  $f \cup h$ 
cannot be monotone. Hence the remaining case is the following.

{\it  Case 2.}  $f \cup h^\prime$ is a monotone function and  $f \cup h$ is a 
function 
but not monotone. --- In this case, there are  $r$,  $u \in T$  such that 
$r <_T u$ 
and  $f(r) = 1$,  $h(u) = 0$. For otherwise, we have   $r <_T u$  satisfying  
$h(r) = 1$,  $f(u) = 0$. It follows that  $u \in w \subseteq Y$,  $r \in Y$  
since  $Y$  is a subtree of  $T$, and  $r \in \text{dom } h^\prime$, 
contradicting the fact that  $f \cup h^\prime$  is monotone.

Now  $r \in w \subseteq Y$  and  $u \in \{s_1, \dots,s_n, t_1, \dots, 
t_m\}$. 
In fact,  $u = t_j$  for some  $j$, since  $u = s_i$  would imply that  
$ u \in \text{dom } h^\prime$, but $f \cup h^\prime$  was monotone. But 
then 
$r \in \text{pred } t_j \cap Y \subseteq Z$,  $r \in \text{dom } f^\prime$, 
and  $f^\prime \cup h$  is not monotone. Thus  
$b \cdot q_{f^\prime} = q_h \cdot q_{f^\prime} = 0$  and  $x = q_f \leq 
q_{f^\prime}$  
show that  $q_{f^\prime}$  separates  $x$  and  $b$.\qed        
\enddemo

\enddocument